\documentclass[a4paper,12pt]{article}
\usepackage[cp1251]{inputenc}
\usepackage[english]{babel}
\usepackage[tbtags]{amsmath}
\usepackage{amsfonts,amssymb,mathrsfs,amscd}

\usepackage{wrapfig}
\newtheorem{theorem}{Theorem}[section]
\newtheorem{lemma}{Lemma}[section]
\newtheorem{definition}{Definition}[section]

 \date{}
 \numberwithin{equation}{section}

\begin{document}

\large
 \centerline{\bf Discrete Riesz MRA  on  local fields of positive characteristic}
  \centerline{\bf G.S.~Berdnikov, S.\,F.~Lukomskii}
%\end{document}
\noindent
 N.G.\ Chernyshevskii Saratov State University\\
 LukomskiiSF@info.sgu.ru\\
 vam21@yandex.ru

\begin{abstract}
  We propose a  method to construct Riesz MRA  on  local fields of positive characteristic and corresponding  scaling step functions that generate it \\
   Bibliography: 20 titles.
\end{abstract}
\noindent
 keywords: local fields, multiresolution
analysis, Riesz wavelet bases, trees.

%---------------------------------------------

\section{Introduction}\label{s1}
  The simplest example  of a local field with positive characteristic is a Vilenkin group. More precisely, Vilenkin group is an additive group of the local field $F^{(s)}$ when $s=1$.  Additive group $F^+$ of the  local field $F^{(s)}$ with positive characteristic is a product $\mathfrak{G}^s$ of Vilenkin groups $\mathfrak{G}$. Therefore discrete wavelets on local fields  is an alternative method for multidimensional discrete data processing.

  V. Protasov and  Yu. Farkov \cite{PF}-\cite{YuF2}  obtained the necessary and sufficient conditions under which a refinable function $\varphi$ generates an orthogonal MRA on Vilenkin group $F^{(1)}$ and indicated some methods for constructing such  refinable functions. They proved that the refinable function $\varphi$  generates an orthogonal
MRA if and only if the mask $m_0$ does not have any blocked sets. The problem of finding the blocked sets is an exhaustive search problem.
 In  articles  \cite{LSBG},\cite{LBK} a new
  method for constructing refinable  step functions is proposed.
  This method is based on a new concept of N-valid trees.
  Apparently, this method gives all step functions generating an orthogonal MRA.
In \cite{YuF3},\cite{FR2} some algorithms for constructing biorthogonal compactly
supported wavelets on Vilenkin groups are researched and new examples of biorthogonal compactly supported wavelets on Vilenkin groups are given. In \cite{SLuk2},\cite{SLuk3}  1-valid  trees are  used for constructing Riesz bases on Vilenkin groups.

The simplest example  of a local field with characteristic zero  is the field $\mathbb Q_p$ of $p$-adic numbers. Wavelet theory over the field $\mathbb Q_p$ is different from the  wavelet theory on Vilenkin groups \cite{SVK}-\cite{AES}.

   In this article we will discuss MRA on local fields of positive characteristic. First results on wavelet analysis on local fields   are received by Chinese mathem\-aticians
  Huikun Jiang, Dengfeng Li, and Ning Jin   \cite{JLJ}. They  introduced the notion of orthogonal MRA on local
  fields, for the fields  $ F^{(s)}$
  of positive characteristic $p$,  proved some properties  and
  gave an algorithm for constructing wavelets for a known scaling
  function. Using these results, they constructed orthogonal MRA and
  corresponding  wavelets for the case when a scaling function is the characteristic  function of a unit ball $\cal D $. Such MRA is  usually called  "Haar MRA" \   and corresponding
  wavelets are called "Haar wavelets".
      In the article \cite{BJ2} Biswaranjan Behera and Qaiser Jahanthe  proved that a function $\varphi\in L^2( F^{(s)})$ is
  a scaling function for MRA in $L^2( F^{(s)})$ if and only if
  \begin{equation}\label{eq1.1}
    \sum_{k\in \mathbb N_0}|\hat\varphi(\xi+u(k))|^2=1\  for\  a.e. \ \xi \in
    {\cal D},
  \end{equation}
  \begin{equation}\label{eq1.2}
    \lim\limits_{j\to\infty}|\hat\varphi(\mathfrak p^j\xi)|=1 \ for\  a.e.\  \xi \in
    F^{(s)},
  \end{equation}
  and there exists an integral periodic function $m_0 \in L^2(\cal
  D)$ such that
  \begin{equation}\label{eq1.3}
  \hat\varphi(\xi)=m_0(\mathfrak p\xi)\hat\varphi(\mathfrak p\xi)\  for \ a.e.\  \xi
  \in F^{(s)}
  \end{equation}
    where $\{u(k) \}$ is the set of shifts, $\mathfrak p$ is a prime element.
  The condition   (\ref{eq1.3}) is the necessary condition for inclusion $V_0\subset V_1$. The condition   (\ref{eq1.2}) is the necessary and sufficient condition for convergence of the product $\prod m_0(\mathfrak p^j\xi)$. The condition   (\ref{eq1.1}) is the necessary and sufficient   condition for orthogonality of the shifts system $\varphi(\xi+u(k))$. It is a difficult problem. If $\hat\varphi(\xi)=1$ on some ball $|\xi|<r$ then we have  only two conditions (\ref{eq1.1}),(\ref{eq1.3}). For such step function
  methods for constructing orthogonal refinable functions are obtained in  the article  \cite{SLAV}. B.Behera and Q.Jahan  \cite{BJ3}
 found a condition on the scaling functions   $\varphi$ and $\tilde\varphi$ for dual MRAs under which
the associated wavelets $\psi_l$ and $\tilde\psi_l$ generate
the biorthogonal affine systems
$(\psi)_{ljk}$ and $(\tilde\psi)_{l,j,k}$ that form Riesz bases for $L_2(F^{(s)})$. Currently, methods for constructing nonorthogonal wavelets on local fields are missing.

In this article we will consider a Riesz MRA with  the step scaling functions $\varphi$ for which  ${\rm supp}\ \hat\varphi$    is obtained by {\it spreading} the unit ball  $K_0^\bot$. We will give an algorithm for constructing such Riesz scaling functions.  To construct this Riesz scaling functions we will use N-valid trees.

\section{Local field of positive  characteristic as a vector space over a finite field}
In works \cite{JLJ}-\cite{BJ3} authors
 use the notation and methods of the book by Taibleson \cite{MT}. We will use another methods \cite{SLAV}.

Let $K=F^{(s)}$ be a local field with positive characteristic $p$.
Its elements are infinite sequences
$$
a=(\dots , {\bf 0}_{n-1}, {\bf a}_n, {\bf a}_{n+1},\dots), {\bf
a}_j\in GF(p^s)
$$
where
$$
 {\bf a}_j=(a_j^{(0)},a_j^{(1)},\dots,a_j^{(s-1)}),\ a_j^{(\nu)}\in
 GF(p).
$$
Let $\lambda\in GF(p^s)$.  Since
$$
 \lambda a=(\dots {\bf 0}_{-1},\lambda,{\bf 0}_1,\dots)\cdot(\dots
{\bf 0}_{n-1},{\bf a}_n,{\bf a}_{n+1},\dots)=
$$
$$
 (\lambda+{\bf 0}x+{\bf 0}x^2+\dots)({\bf a}_nx^n+{\bf a}_{n+1}x^{n+1}+\dots)=\lambda {\bf a}_nx^n+\lambda
 {\bf a}_{n+1}x^{n+1}+\dots=
 $$
 $$
 =(\dots {\bf 0}_{n-1},\lambda {\bf a}_n,\lambda {\bf a}_{n+1},\dots)
$$
it follows that the product  $\lambda  a$ is defined
coordinate-wise.  The sum in the field $K$ is defined
coordinate-wise also.  With such operations   $F^{(s)}$ is a vector
space. If we define the modulus $|\lambda|$ by the equation
$$
 |\lambda|=\left\{
           \begin{array}{ll}
 1,& \lambda \ne 0,\\
 0,& \lambda = 0.\\
  \end{array} \right.
$$
and norm $\|a\|$ by the equation
\begin{equation}\label{eq3.2}
 \|a\|=\frac{1}{p^{sn}},\;{\bf a}_n\ne 0
\end{equation}
then we can consider the field  $F^{(s)}$ as a vector normalized
space over the field $GF(p^s)$.
Let $$
K_n=\{a=(\dots , {\bf 0}_{n-1}, {\bf a}_n, {\bf a}_{n+1},\dots):n\in \mathbb N, {\bf
a}_j\in GF(p^s)\}
$$
be a ball of radius $\frac{1}{p^{sn}}$. For any $n\in\mathbb Z$ choose
an element $g_n\in K_n\setminus K_{n+1} $ and fix it. We will call this system $(g_n)_{n\in \mathbb Z }$  a basic sequence.
\begin{theorem}[\cite{SLAV}]
Let $(g_n)_{n\in \mathbb Z}$ be a fixed basic sequence in $K$. Any element $a\in K$ may by
written as sum of the series
\begin{equation}\label{eq3.3}
a=\sum_{n\in \mathbb Z} \overline{\lambda}_ng_n, \
\overline{\lambda}_n\in GF(p^s).
\end{equation}
\end{theorem}
It means that the sequence $(g_n)$ is a basis of vector space $K$. Further we will suppose $g_n=(\dots , {\bf 0}_{n-1}, {\bf 1}_n, {\bf 0}_{n+1},\dots)$, where ${\bf 1}_n=(1,0,...,0)$.

\begin{definition}
 The operator
 $$
 \mathcal{A}:a=\sum_{n\in \mathbb Z}\overline{\lambda}_ng_n\longmapsto \sum_{n\in \mathbb
 Z}\overline{\lambda}_ng_{n-1}
 $$
 is called a delation operator.
\end{definition}
 {\bf Remark 1.}  Since additive group
 $K^+=F^{(s)+}$ is Vilenkin group $\mathfrak G$ with $\mathfrak G_{ns}=F_n^{(s)+}$
 it follows that $\mathcal{A}K_n=\mathcal{A}K_{n-1} $ and
 $\int\limits_{K^+}f(\mathcal{A}u)\,d\mu=\frac{1}{p^s}\int\limits_{K^+}f(x)\,d\mu$.

\section{Set of characters as vector space over a finite field}
Now we define Rademacher function on the vector space  $K$.
 If
$$
a=(\dots , {\bf 0}_{k-1}, {\bf a}_k, {\bf a}_{k+1},\dots), {\bf
a}_j\in GF(p^s)
$$
and
$$
 {\bf a}_j=(a_j^{(0)},a_j^{(1)},\dots,a_j^{(s-1)}),\ a_j^{(\nu)}\in
 GF(p)
$$
 then we define functions  $r_n(a)=e^{\frac{2\pi i}{p}a_k^{(l)}}$, where $n=ks+l$, $0\leq l<s$.
 \begin{lemma}[\cite{SLAV}]
Any character $\chi\in X$ can be expressed uniquely as product
\begin{equation}\label{eq3.1}
  \chi=\prod_{n=-\infty}^{+\infty}
  r_n^{\alpha_n}\;\;(\alpha_n=\overline{0,p-1}),
   \end{equation}
 in which the number of factors with positive numbers are finite.
 \end{lemma}
   If we write the character $\chi$ as
  $$
  \chi =\prod_{k\in \mathbb Z}r_{ks+0}^{a_k^{(0)}}
  r_{ks+1}^{a_k^{(1)}}\dots  r_{ks+s-1}^{a_k^{(s-1)}}
  $$
  and denote
  $$
  {\bf r}_k^{{\bf a}_k}:=r_{ks+0}^{a_k^{(0)}}
  r_{ks+1}^{a_k^{(1)}}\dots r_{ks+s-1}^{a_k^{(s-1)}},
  $$
  where ${\bf a}_k=(a_k^{(0)},a_k^{(1)},...,a_k^{(s-1)})\in
  GF(p^s)$ then we  can write the character $\chi$ as the product
\begin{equation}\label{eq3.2}
     \chi =\prod_{k\in \mathbb Z}{\bf r}_k^{{\bf a}_k}.
    \end{equation}
  The function ${\bf r}_k={\bf r}_k^{(1,0,\dots,0)}$ is
  called Rademacher function.

 Assume by definition
  $$
   ({\bf r}_k^{{\bf a}_k})^{{\bf b}_k}:={\bf r}_k^{{\bf a}_k{{\bf b}_k}},\ {\bf a}_k, {\bf b}_k\in
  GF(p^s).
   $$
 and
    $$
    \chi^{\bf b}:=\prod\limits_{k\in\mathbb Z} ({\bf r}_k^{{\bf a}_k})^{\bf
    b}.
    $$
 Then
 $$
 {\bf r}_k^{{\bf u}\dot+{\bf v}}={\bf r}_k^{{\bf u}}\cdot{\bf r}_k^{\bf
 v},\;{\bf u},{\bf v}\in GF(p^s)
 $$
  and the set of characters of the field $F^{(s)}$ is a vector space
   $(X,\; *,\; \cdot^{GF(p^s)})$ over  the finite field  $GF(p^s)$ with product as interior
   operation and power as exterior operation.
      It follow from \eqref{eq3.2} that annihilator  $(F_k^{(s)})^\bot$ consists
   of characters $\chi={\bf r}_{k-1}^{{\bf a}_{k-1}}{\bf r}_{k-2}^{{\bf
   a}_{k-2}}...$. \\
  The next lemma is the basic property  of Rademacher
  functions    on local field with positive characteristic.
   \begin{lemma}[\cite{SLAV}] Let
   $g_j=(\dots,{\bf 0}_{j-1},(1,0,\dots,0)_j,{\bf
   0}_{j+1},\dots)\in F^{(s)}$, ${\bf a}_k,{\bf u}\in GF(p^s)$. Then
   $({\bf r}_k^{{\bf a}_k},{\bf u}g_j)=1$ for any $k\ne j$.
   \end{lemma}
   \section{Riesz MRA on local fields of positive characteristic}

     Denote
     $$
     H_0
  =\{a: a={\bf a}_{-1}g_{-1}\dot+\dots\dot+{\bf
  a}_{-\nu}g_{-\nu},\  \nu\in \mathbb N,\ {\bf a}_{-j}\in GF(p^s)\},
 $$
   $$
     H_0^{(\nu)}
  =\{a: a={\bf a}_{-1}g_{-1}\dot+\dots\dot+{\bf
  a}_{-\nu}g_{-\nu},\  {\bf a}_{-j}\in GF(p^s)\},\  \nu\in \mathbb N.
 $$
     $H_0$ is an analog of the set $(N_0)^s=N_0\times\dots\times N_0$ .

   Define a dilation operator $\mathcal{A}$ on the set of characters by
    the equation $(\chi \mathcal{A},x)=(\chi, \mathcal{A}x)$.
 It is evident  that ${\bf r}_j\mathcal{A}={\bf r}_{j+1}, \
    (K_n^+)^\bot\mathcal{A}=(K_{n+1}^+)^\bot $ and
 $\int\limits_{X}f(\chi\mathcal{A})\,d\nu=
 \frac{1}{p^s}\int\limits_Xf(\chi)\,d\nu$ \cite{SLAV}.

  We will use next properties of annihilators \cite{SLAV}.\\
    1) $\int\limits_{(K_n^+)^\bot}(\chi,x)\,d\nu(\chi)=p^{sn}{\bf 1}_{K_n^+}(x)$,\\
  2) $\int\limits_{K_n^+}(\chi,x)\,d\mu(x)=\frac{1}{p^{sn}}{\bf 1}_{(K_n^+)^\bot}(\chi)$.\\
   3) Let  $\chi_{n,l}={\bf r}_n^{{\bf a}_n}{\bf r}_{n+1}^{{\bf
 a}_{n+1}}\dots {\bf r}_{n+l}^{{\bf a}_{n+l}}$ be a character which does
 not belong to $({K_n^+})^\bot$. Then
 $$
 \int\limits_{({K_n^+})^\bot\chi_{n,l}}(\chi,x)\,d\nu(\chi)=p^{ns}(\chi_{n,l},x){\bf
 1}_{K_n^+}(x).
 $$
 4)
 Let $h_{n,l}={\bf
 a}_{n-1}g_{n-1}\dot+{\bf a}_{n-2}g_{n-2}\dot+\dots\dot+{\bf a}_{n-l}g_{n-l}\notin
 K_n^+$. Then
 $$
 \int\limits_{K_n^+\dot+h_{n,l}}(\chi,x)\,d\mu(x)=\frac{1}{p^{ns}}(\chi,h_{n,l}){\bf
 1}_{({K_n^+})^\bot}(\chi).
 $$

\begin{definition}
  A family of closed subspaces $V_n$, $n\in\mathbb Z$,
 is said to be a~Riesz multi\-resolution analysis  of~$L_2(K)$
 if the following  axioms (conditions) are satisfied:
 \begin{itemize}
 \item[A1)] $V_n\subset V_{n+1}$;
 \item[A2)] ${\vrule width0pt
 depth0pt height11pt} \overline{\bigcup_{n\in\mathbb
 Z}V_n}=L_2(K)$ and $\bigcap_{n\in\mathbb Z}V_n=\{0\}$;
 \item[A3)] $f(x)\in V_n$  $\Longleftrightarrow$ \ $f({\cal A} x)\in V_{n+1}$ (${\cal A}$~is a~dilation
 operator);
 \item[A4)] $f(x)\in V_0$ \ $\Longrightarrow$ \
 $f(x\dot - h)\in V_0$ for all $h\in H_0$;
  \item[A5)] there exists
 a~function $\varphi\in L_2(K)$ such that the system
 $(\varphi(x\dot - h))_{h\in H_0}$ is a Riesz basis
 for~$V_0$.
\end{itemize}
 A function~$\varphi$ occurring in axiom~A5 is called
a~\textit{scaling function}.
\end{definition}
We recall  that the family $(f_j)\subset L_2(K)$ is called a Riesz system with constants $A$ and $B\  (A,B>0)$ if, for every sequence $C=(c_j)\in l_2$ the series
$\sum_j c_jf_j$ converges in $L_2(K)$ and
\begin{equation}                                \label{eq4.1}
 A||C||_{l_2}^2\le \left\|\sum
c_n\varphi_n\right\|_{L_2(G)}^2\le B||C||_{l_2}^2.
\end{equation}
% Note that condition (\ref{eq4.1}) already implies that $\sum_j c_jf_j$ %converges. It is also clear that every  Riesz system is a basic in the space
%$$
%V:=\left\{f=\sum c_n\varphi_n:\;\sum|c_k|^2<+\infty\right\}.
%$$

\begin{lemma}[ \cite{BJ2}, {\rm Theorem 4.1.}]
 Let $(V_n)_{n\in \mathbb Z}$  be a sequence of closed subspaces of $L_2(K)$ satisfying
conditions (A1), (A3) and (A5) of Definition 4.1. Then, $\bigcap\limits_{n\in \mathbb Z}V_n=\{0\}$.
\end{lemma}

\begin{lemma}[ \cite{BJ2}, {\rm Theorem 4.2.}]Let $\hat\varphi\in L_2(X)$ be a continuos function at the point $\chi=1$ and $\hat\varphi(1)\ne 0$.
 Suppose also that $(V_n)_{n\in \mathbb Z}$  is a sequence of closed subspaces of $L_2(K)$ satisfying
conditions (A1), (A3) and (A5) of Definition 4.1. Then
$\overline{\bigcup\limits_{n\in \mathbb Z}V_n}=L_2(K)$.
\end{lemma}

Next we will follow the conventional approach. Let
 $\varphi(x)\,{\in}\, L_2(K)$, and suppose that
 $(\varphi(x\dot -\nobreak h))_{h\in H_0}$ form a Riesz basis in the closure of their linear hull  in the norm $L_2(K)$. With the function~$\varphi$ and the
 dilation operator~${\cal A}$, we define
 subspaces $V_n=\overline{{\rm span } (\varphi({\cal A}^nx\dot - h))_{h\in H_0}}$ closed in $L_2(K)$ .  If the family $(V_n)_{n\in \mathbb Z}$ is an MRA, we will say that the function $\varphi$ generates MRA. It is clear that the
  system $(p^{\frac{ns}{2}}\varphi({\cal A}^n x\dot-h))_{h\in H_0}$ is
 a Riesz basis for $ V_n$ and $f(x)\in V_n$ if and only if $f({\cal A}x)\in V_{n+1}$. We want to propose an algorithm for constructing a function $\varphi$ that generate  a Riesz MRA and corresponding wavelets.

 We shall assume that the function $\varphi$ generating a Riesz MRA satisfies the inequality
 \begin{equation}                                \label{eq4.2}
0< A\le |\hat\varphi(\chi)|^2\le B,
\end{equation}
 on some set $E\subset X$ of measure $\nu E =1$ that is obtained by {\it spreading} the set $K_0^\bot$. We now give a precise characterization or $E$.
 \begin{definition}
Let  $K=F^{(s)}$ be a local field of characteristic $p$, $ X$--group of additive characters for $K^+$, $N\in\mathbb N, M\in\mathbb N_0=\mathbb N\bigsqcup\{0\}$. A set  $E\subset X$
is said to be
$(N,M)$-elementary if it is the disjoint union of $p^{Ns}$ cosets of the form
$$
 (K_{-N}^+)^\bot\zeta_j=
(K_{-N}^+)^\bot\underbrace{{\bf r}_{-N}^{\overline{\alpha}_{-N}}{\bf r}_{-N+1}^{\overline{\alpha}_{-N+1}}\dots
{\bf r}_{-1}^{\overline{\alpha}_{-1}}}_{\xi_j}
\underbrace{{\bf r}_{0}^{\overline{\alpha}_{0}}{\bf r}_{1}^{\overline{\alpha}_{1}}\dots
{\bf r}_{M-1}^{\overline{\alpha}_{M-1}}}_{\eta_j}=
(K_{-N}^+)^\bot\xi_j\eta_j,
$$
for $j=0,1,\dots,p^{Ns}-1$ such that the following conditions hold\\
1) $\bigsqcup\limits_{j=0}^{p^{Ns}-1}
(K_{-N}^+)^\bot\xi_j=(K_{0}^+)^\bot$,
 $(K_{-N}^+)^\bot\xi_0=(K_{-N}^+)^\bot$.\\
2) For every $l=\overline{0,M+N-1}$ we have $((K_{-N+l+1}^+)^\bot\setminus (K_{-N+l}^+)^\bot)\bigcap
E\ne\emptyset$.
\end{definition}
So, to obtain an $(N,M)$-elementary set, we shift any coset\\ $(K_{-N}^+)^\bot
{\bf r}_{-N}^{\overline{\alpha}_{-N}}{\bf r}_{-N+1}^{\overline{\alpha}_{-N+1}}\dots
{\bf r}_{-1}^{\overline{\alpha}_{-1}}$ per the unique element
${\bf r}_{0}^{\overline{\alpha}_{0}}{\bf r}_{1}^{\overline{\alpha}_{1}}\dots
{\bf r}_{M-1}^{\overline{\alpha}_{M-1}}$
that any difference $((K_{-N+l+1}^+)^\bot\setminus (K_{-N+l}^+)^\bot)$ contains at least one shift.

 \begin{lemma}[\cite{SLAV}]
 The set $H_0\subset  K$ is a total orthonormal  system on any
 $(N,M)$-elementary set $E\subset X$.
 \end{lemma}
 \begin{lemma}
 Let $K=F^{(s)}$ be a local field of characteristic $p$,
 $E\subset (K_M^+)^\bot$ an $(N,M)$-elementary set in
 $X$, $\varphi\in L^2(K)$, ${\rm supp}\ \hat\varphi =E$,  $A,B>0$.
 The system of shifts  $(\varphi(x\dot- h))_{h\in H_0}$ is a Riesz system with constants $A$ and $B$ if and only if
 $$
 A\le |\hat\varphi(\chi)|^2\le B,
 $$
  a.e. on $E$.
 \end{lemma}

 {\bf Proof.} S~u~f~f~i~c~i~e~n~c~y. First, we find an upper bound for $\left\|\sum\limits_{h\in\tilde
 H_0}c_h\varphi(x\dot- h)\right\|_2^2$ assuming that
 $\tilde H_0\subset H_0$  is a finite set. Using Plancherel's equality, we have
 $$
 \int\limits_{K}\left|\sum\limits_{h\in\tilde
 H_0}c_h\varphi(x\dot- h)\right|^2\,d\mu(x)=\int\limits_{X}\left|\sum\limits_{h\in
 \tilde H_0}c_h\hat\varphi(\chi)\overline{(\chi,h)}\right|^2\,d\nu(\chi)\le
 $$
 $$
 B\int\limits_{E}\left|\sum\limits_{h\in \tilde
 H_0}c_h\overline{(\chi,h)}\right|^2\,d\nu(\chi)=
 B\sum_{j=0}^{p^{Ns}-1}\int\limits_{(K_{-N}^+)^\bot\zeta_j}\left|\sum\limits_{h\in\tilde
 H_0}c_h\overline{(\chi,h)}\right|^2\,d\nu(\chi).
 $$
 We rewrite the inner integral using invariance of the integral

 $$
 \int\limits_{(K_{-N}^+)^\bot\zeta_j}\left|\sum\limits_{h\in\tilde
 H_0}c_h\overline{(\chi,h)}\right|^2\,d\nu(\chi)=\int\limits_{X}{\bf
 1}_{(K_{-N}^+)^\bot\zeta_j}(\chi)\left|\sum\limits_{h\in\tilde
 H_0}c_h\overline{(\chi,h)}\right|^2\,d\nu(\chi)=
 $$
 $$
 =\int\limits_{X}{\bf 1}_{(K_{-N}^+)^\bot\zeta_j}(\chi\eta_j)\left|\sum\limits_{h\in\tilde
 H_0}c_h\overline{(\chi\eta_j,h)}\right|^2\,d\nu(\chi)=
 $$
 $$
 =\int\limits_{X}{\bf 1}_{(K_{-N}^+)^\bot\xi_j}(\chi)\left|\sum\limits_{h\in\tilde
 H_0}c_h\overline{(\chi\eta_j,h)}\right|^2\,d\nu(\chi)=
 $$
 $$
 =\int\limits_{(K_{-N}^+)^\bot
 \xi_j}\sum\limits_{h\in\tilde
 H_0}\sum\limits_{g\in\tilde
 H_0}c_h\overline{c_g}\overline{(\eta_j,h)}(\eta_j,g)\overline{(\chi,h)}(\chi,g)\,d\nu(\chi).
 $$

 Since
 $$
 (\eta_j,h)=({\bf r}_0^{\overline{\alpha}_0}{\bf r}_1^{\overline{\alpha}_1}\dots
 {\bf r}_{M-1}^{\overline{\alpha}_{M-1}},a_{-1}{g}_{-1}\dot+a_{-2}{ g}_{-2}\dot+\dots\dot+a_{-\nu}{ g}_{-\nu})=1,
 $$
 $$
 (\eta_j,g)=({\bf r}_0^{\overline{\alpha}_0}{\bf r}_1^{\overline{\alpha}_1}\dots
 {\bf r}_{M-1}^{\overline{\alpha}_{M-1}},b_{-1}{ g}_{-1}\dot+b_{-2}g_{-2}\dot+\dots\dot+b_{-\nu}g_{-\nu})=1,
 $$
 we get
 $$
 \sum_{j=0}^{p^N-1}\int\limits_{(K_{-N}^+)^\bot\zeta_j}\left|\sum\limits_{h\in\tilde
 H_0}c_h\overline{(\chi,h)}\right|^2\,d\nu(\chi)=\sum_{g,h\in\tilde
 H_0}c_h\overline{c_g}\sum_{j=0}^{p^N-1}\int\limits_{(K_{-N}^+)^\bot\xi_j}\overline{(\chi,h)}(\chi,g)\,d\nu(\chi)=
 $$
 $$
 =\sum_{g,h\in\tilde H_0}c_h\overline{c_g}\int\limits_{(K_{0}^+)^\bot}\overline{(\chi,h)}(\chi,g)\,d\nu(\chi)=\sum_{h\in\tilde
 H_0}|c_h|^2,
 $$
 and, therefore,
 $$
 \left\|\sum_{h\in\tilde H_0}c_h\varphi(x\dot + h)\right\|_2^2\le
 B\sum_{h\in\tilde H_0}|c_h|^2.
 $$

 Similarly we obtain that
 $$
 \left\|\sum_{h\in\tilde H_0}c_h\varphi(x\dot +h)\right\|_2^2\ge
 A\sum_{h\in\tilde H_0}|c_h|^2.
 $$
 N~e~c~e~s~s~i~t~y. Let $(\varphi(x\dot-h))_{h\in H_0}$ be a Riesz system with bounds A and B i.e. for
 any $(c_h)\in l_2$ the inequality
 \begin{equation}    \label{eq4.3}
 A\sum_{h\in H_0}|c_h|^2\le \left\| \sum_{h\in H_0}c_h\varphi(x\dot-h)\right\|^2_{L_2(K)}\le B\sum_{h\in H_0}|c_h|^2
\end{equation}
holds. It follows that
$$
f(x)=\sum_{h\in H_0}c_h\varphi(x\dot-h)\in L_2(K).
$$
By Plancherel's equality we have
$$
\left\|\sum_{h\in H_0}|c_h\varphi(x\dot-h)\right\|_{L_2(K)}^2=\left\|\sum_{h\in H_0}|c_h\hat\varphi(\chi)\overline{(\chi,h)}\right\|_{L_2(E)}^2.
$$
Therefore we can write the equality (\ref{eq4.3})
in the form
$$
A\sum_{h\in H_0}|c_h|^2\le \hat\varphi(\chi)\| \sum_{h\in H_0}c_h\overline{(\chi,h)}\|^2_{L_2(E)}\le B\sum_{h\in H_0}|c_h|^2.
$$
Denote $g(\chi)=\sum\limits_{h\in H_0}c_h\overline{(\chi,h)}$. Since $(H_0)$ is an orthonormal system in $L_2(E)$ we obtain
$$
A\|g\|_{L_2(E)}^2\le \int\limits_E|\hat\varphi(\chi)|^2|g(\chi)|^2d\,\nu(\chi)\le B\|g\|_{L_2(E)}^2
$$
or another
$$
A\le \int\limits_E|\hat\varphi(\chi)|^2\frac{|g(\chi)|^2}{\|g\|^2}d\,\nu(\chi)\le B.
$$
It follows that for any $h\in L_2(E)$ with norm $\|h\|_{L_2(E)}=1$
\begin{equation}                            \label{eq4.4}
%\label{eleventh4.3}
A\le \int\limits_E|\hat\varphi(\chi)|^2|h(\chi)|d\,\nu(\chi)\le B.
\end{equation}
Therefore $A\le {\rm vrai~ sup} |\hat\varphi(\chi)|^2\le B $.

Let us show that $A\le |\hat\varphi(\chi)|^2$ a.e. in $E$. Assume the converse. Then there exists $E_1\subset E$ such that $\nu(E_1)>0$
and $|\hat\varphi(\chi)|^2<A-\varepsilon$ on the set $E_1$. Taking $h(x)=\frac{{\bf 1}_{E_1}(\chi)}{\nu(E_1)}$ we have
$$
\int\limits_E|\hat\varphi(\chi)|^2|h(\chi)|d\,\nu(\chi)=\int\limits_{E_1}|\hat\varphi(\chi)|^2\frac{1}{\nu(E_1)}d\,\nu(\chi)\le A-\varepsilon<A.
$$
But this contradicts inequality \ref{eq4.4}. $\square$

\begin{lemma}
Let $K=F^{(s)}$ be a local field, $E\subset X$ an $(N,M)$ elementary set in  $X$, ${\rm
supp}\,\hat\varphi={\rm supp}\,\hat\psi=E$.\\
1) If $\hat\varphi(\chi)\overline{\hat\psi(\chi)}=1$ a.e. on
$E$, then  $(\varphi(x\dot-h),\psi(x\dot-g))_{h,g\in H_0}$ is a biorthonormal system on $E$.\\
2) Conversely, if $(\varphi(x\dot-h),\psi(x\dot-g))_{h,g \in H_0}$ is
a biorthonormal system on $E$, then
$\hat\varphi(\chi)\overline{\hat\psi(\chi)}=1$ a.e. on $E$.
\end{lemma}
{\it Proof.} 1) Using Plancherel’s equality and Lemma 4.1, we have
$$
\int\limits_{K}\varphi(x\dot-h)\overline{\psi(x\dot-g)}\,d\mu(x)=\int\limits_E\overline{(\chi,h)}(\chi,g)\,d\nu(\chi)=\delta_{h,g}.
$$
2) Using Plancherel’s equality and  a biorthogonality of the system\\
$(\varphi(x\dot-h),\psi(x\dot-h))$ we have
$$
\int\limits_{K}\overline{\varphi(x\dot-h)}\psi(x)\,d\mu(x)=\int\limits_E\overline{\hat\varphi(\chi)}\hat\psi(\chi)(\chi,h)\,d\nu(\chi)=
c_h=\left\{\begin{array}{ll}
  1, & h=0 \\
  0 & h\ne0,h\in H_0 \\
\end{array}
\right.
$$
By the uniqueness theorem
$\overline{\hat\varphi(\chi)}\hat\psi(\chi)=1$ a.e. on $E$.
$\square$\\
The following lemma obviously follows from the equality
$\int\limits_Kf({\cal
A}x)d\,\mu=\frac{1}{p^s}\int\limits_Kf(x)d\,\mu$.\\

\begin{lemma}
 Let $n\in \mathbb N$. The shift system $(\varphi(x\dot-h))_{h\in H_0}$ is
 a Riesz system with constants $A$ and $B$ if and only if the system $(p^{\frac{ns}{2}}\varphi({\cal A}^n x\dot-h))_{h\in H_0}$ is
 a Riesz system with constants $A$ and $B$.
\end{lemma}
Recall that  $V_n=\overline{{\rm span } (\varphi({\cal A}^nx\dot - h))_{h\in H_0}}$.
 \begin{lemma}
  Suppose that $\varphi\in L_2(K)$, $E$ is an $(N,M)$-elementary set, ${\rm supp}\ \hat \varphi = E$,
 $\hat \varphi$ satisfies the conditions
 (\ref{eq4.2}) on  $E$. Then  $V_0\subset V_1$ if and only if the  function $\varphi$
  satisfies the equation
   \begin{equation}                                \label{eq4.5}
\varphi(x)=\sum_{h\in H_0}\beta_h\varphi({\cal
A}x\dot-h),\;\;\sum_{h\in H_0}|\beta_h|^2<+\infty.
 \end{equation}
  \end{lemma}
 {\it Proof.}  N~e~c~e~s~s~i~t~y. By lemma 4.4 $\varphi(x\dot-h)$ is
 a Riesz system with constants $A$ and $B$. Taking  Lemma 4.6 into account, we get that   $\sqrt{p^s}(\varphi({\cal A}x\dot-h))_{h\in H_0}$ is a
Riesz system with constants  $A$ and $B$, and, therefore, form a basis of
$V_1$. Since $\varphi\in V_0\subset V_1$, the equation (\ref{eq4.5}) holds and we have $\sum|\beta_h|^2<+\infty$. \\
 S~u~f~f~i~c~i~e~n~c~y.
 Let $f\in{\rm span}\,(\varphi(x\dot-h))_{h\in
 H_0}$, i.e.
\begin{equation}                                \label{eq4.6}
f(x)=\sum_{\tilde h\in\tilde H\subset H_0}\alpha_{\tilde
h}\varphi(x\dot-\tilde h),
\end{equation}
  where $\tilde H\subset H_0$ is the finite set. Substituting  (\ref{eq4.5})
in (\ref{eq4.6}), we obtain
\begin{equation}                                \label{eq4.7}
f(x)=\sum_{\tilde h\in\tilde H}\alpha_{\tilde h}\sum_{h\in
H_0}\beta_h\varphi({\cal A}x\dot-({\cal A}\tilde
h\dot+h))=\sum_{h\in H_0}\sum_{\tilde h\in\tilde H}\alpha_{\tilde
h}\beta_h\varphi({\cal A}x\dot-({\cal A}\tilde h\dot+h)).
\end{equation}
Since the set $H_0$ is a group and ${\cal A}\tilde h\in H_0$, it follows that $f(x)\in V_1$
$\square$\\
Therefore we need  to look for a~function
 $\varphi\in L_2(K)$, that generates an~MRA
 in~$L_2(K)$, as a~solution of the refinement
 equation (\ref{eq4.5}). A solution of the  refinement equation (\ref{eq4.5}) is called a {\it refinable function}.

\begin{theorem}
  Suppose that $\varphi\in L_2(K)$, $E\subset X$ is an $(N,M)$-elementary set, ${\rm supp}\ \hat \varphi = E$,
 $\hat \varphi$ satisfies the conditions
 (\ref{eq4.2}) on  $E$
. Then  $V_n\subset V_{n+1}$ if and only if the  function $\varphi$
  satisfies the equation (\ref{eq4.5})
     \end{theorem}
This theorem  follows from lemmas 4.6-4.7.
 \begin{theorem}Let $\hat\varphi\in L_2(X)$ be a continuous function at the point $\chi=1$ and $\hat\varphi(1)\ne 0$. Suppose that $E\subset X$ is an $(N,M)$-elementary set, ${\rm supp}\ \hat \varphi = E$,
 $\hat \varphi$ satisfies the conditions
 (\ref{eq4.2}) and (\ref{eq4.5}).
  Then  $\varphi$ generates an Riesz MRA.
  \end{theorem}
{\it Proof.} The property A5) follows from  lemma 4.4. The property A4) is true, as the set $H_0$ is a group. The property A3) is evident. The property A1) follows from  theorem 4.1. The property A2) follows from  lemmas  4.1. and 4.2. $\square$

\section{Construction of scaling function}
In this section we will construct functions $\hat \varphi$ for which the conditions of theorem 4.2 are satisfied.
 The refinement equation (\ref{eq4.5}) may be written in the form
 \begin{equation}                                                        \label{eq5.1}
 \hat\varphi(\chi)=m_0(\chi)\hat\varphi(\chi{\cal
  A}^{-1}),
 \end{equation}
  where

 \begin{equation}                                                        \label{eq5.2}
 m_0(\chi)=\frac{1}{p^s}\sum_{h\in
 H_0}\beta_h\overline{(\chi{\cal A}^{-1},h)}
 \end{equation}
 is a mask of the equation (\ref{eq4.4}).
 We will use equation (\ref{eq5.1}) to construct the refinable function $\varphi$. First we will construct the support of $\hat\varphi$ using a concept of $N$-valid tree \cite{LSBG}.
  \begin{definition}
     Let $T$ be a tree directed from the root on the set of nodes $GF(p^s)$. The tree $T$ is called as N-valid if the following properties are valid:\\
  a)The nodes of this tree are elements $\bar{\alpha}\in GF(p^s)$\\
  b)The root of $T$ is ${\bf 0}=(0,0,...,0)$\\
  c)For any $j=0,1,...,N-1$ the set $T_j$ of  level $j$ nodes is the set $\{{\bf 0}\}$\\
  d)Any path $(\bar \alpha_k\rightarrow \bar \alpha_{k+1}\rightarrow\dots\rightarrow\bar \alpha_{k+N-1})$ of length $N-1$ is present in the tree $T$ exactly one time.
  \end{definition}
   For example for $p=s=N=2$ we can construct the tree

   \begin{picture}(210,90)
  \put(84,15){$(0,0)$}
  \put(90,22){\vector(0,1){10}}

    \put(84,35){$(0,0)$}
   \put(90,42){\vector(0,1){12}}

    \put(91,40){\vector(4,1){60}}
    \put(89,40){\vector(-4,1){60}}
    \put(147,56){$(1,1)$}
     \put(86,56){$(1,0)$}
      \put(24,56){$(0,1)$}

    \put(25,61){\vector(-3,1){30}}
    \put(27,61){\vector(-1,1){10}}
     \put(29,61){\vector(1,1){10}}
      \put(30,61){\vector(2,1){20}}
    \put(-10,74){$(0,0)$}
    \put(10,74){$(0,1)$}
    \put(28,74){$(1,0)$}
    \put(44,74){$(1,1)$}

    \put(88,61){\vector(-3,1){30}}
    \put(90,61){\vector(-1,1){10}}
     \put(92,61){\vector(1,1){10}}
      \put(93,61){\vector(2,1){20}}
    \put(57,74){$(0,0)$}
    \put(75,74){$(0,1)$}
    \put(92,74){$(1,0)$}
    \put(107,74){$(1,1)$}

   \put(152,61){\vector(-3,1){30}}
    \put(154,61){\vector(-1,1){10}}
     \put(156,61){\vector(1,1){10}}
      \put(157,61){\vector(2,1){20}}
    \put(120,74){$(0,0)$}
    \put(138,74){$(0,1)$}
    \put(157,74){$(1,0)$}
    \put(172,74){$(1,1)$}
    \put(70,0){Figure 1}
   \end{picture}

Here  we give a method for construction of N-valid trees for any $N,s,p$. First we construct a {\it basic tree} $T_B$ of smallest height.
 Let  $\bar\alpha_0,\bar\alpha_1,\dots,\bar\alpha_{p^s-1}$ be all elements of the finite field $GF(p^s)$, and $\bar\alpha_0=(0,0,...,0)$. We construct a basic tree $T_B$ in the following way.\\
 1. Choose a path $(\bar \alpha_0\rightarrow \bar \alpha_{0}\rightarrow\dots\rightarrow\bar \alpha_{0})$ of length $N-1$. This path  contains N nodes $\bar \alpha_{0}$ and  the level of the last node is $N-1$ .\\
 2.  Connect all elements $\bar \alpha_{1},\bar \alpha_{2},\dots,\bar \alpha_{p^s-1}$ to the last node . We get a tree of height $H=N$ in which any path of length $N-1$ is present not more than once. But not all paths of length $N-1$ are present in this tree. This tree is shown in Fig. 2.\\

  \begin{picture}(120,50)
  \put(0,30){$\alpha_0$}
  \put(8,30){\vector(1,0){10}}
  \put(20,30){$\alpha_0$}
  \put(28,30){\vector(1,0){10}}
  \put(39,30){$\dots$}
  \put(48,30){\vector(1,0){10}}
  \put(60,30){$\alpha_0$}

\put(68,30){\vector(1,0){10}}
\put(68,32){\vector(1,1){10}}
\put(68,28){\vector(1,-1){10}}
\put(78,42){$\alpha_1$}
\put(80,29){$\alpha_2$}
\put(77,15){$\alpha_{p^s-1}$}

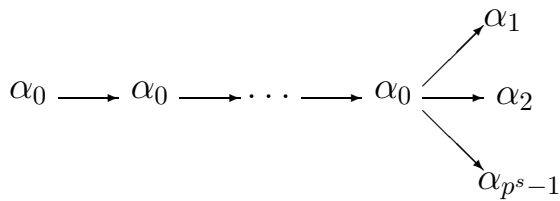
\put(40,0){Figure 2. Tree after 2 steps }
      \end{picture}

    \vspace{0.5cm}
 \noindent
 3.  Now we connect all elements $\bar \alpha_{0},\bar \alpha_{1},\dots,\bar \alpha_{p^s-1}$ to every  node $(\bar \alpha_{1},\bar \alpha_{2},\dots,\bar \alpha_{p^s-1})$ of level $N$ and get a tree of height $H=N+1$. We can see this tree on Fig.3. After $(N+1)$-th step we obtain the $N$-valid three $T_B$ of smallest hight.

  \begin{picture}(120,110)
  \put(0,60){$\alpha_0$}
  \put(8,60){\vector(1,0){10}}
\put(20,60){$\alpha_0$}
  \put(28,60){\vector(1,0){10}}
  \put(39,60){$\dots$}
  \put(48,60){\vector(1,0){10}}
  \put(60,60){$\alpha_0$}

\put(68,60){\vector(1,0){10}}
\put(68,62){\vector(1,2){10}}
\put(68,58){\vector(1,-2){10}}
\put(78,82){$\alpha_1$}
\put(80,59){$\alpha_2$}
\put(67,35){$\alpha_{p^s-1}$}

\put(82,86){\vector(0,1){12}}
\put(85,85){\vector(1,1){12}}
\put(85,82){\vector(1,0){12}}
\put(78,100){$\alpha_0$}
\put(95,100){$\alpha_1$}
\put(98,82){$\alpha_{p^s-1}$}

\put(85,63){\vector(1,1){12}}
\put(87,61){\vector(1,0){10}}
\put(85,59){\vector(1,-1){12}}
\put(98,74){$\alpha_0$}
\put(98,60){$\alpha_1$}
\put(98,45){$\alpha_{p^s-1}$}

\put(81,36){\vector(1,0){14}}
\put(83,34){\vector(1,-1){12}}
\put(81,32){\vector(0,-1){10}}
\put(96,35){$\alpha_0$}
\put(96,20){$\alpha_1$}
\put(78,18){$\alpha_{p^s-1}$}
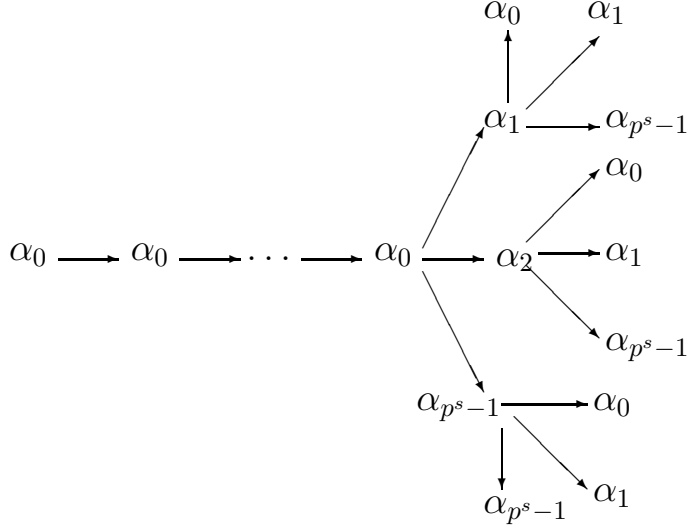
\put(40,5){Figure 3. Tree after 3 steps }
      \end{picture}
 \begin{picture}(40,40)
  \put(20,50){$\alpha_1$}
 \put(28,50){\vector(1,0){10}}
\put(28,52){\vector(1,1){10}}
\put(28,48){\vector(1,-1){10}}
\put(38,62){$\alpha_0$}
\put(40,49){$\alpha_1$}
\put(37,35){$\alpha_{p^s-1}$}
\put(20,5){Figure 4. Subtree}
\end{picture}

  \begin{picture}(120,100)
  \put(0,60){$\alpha_0$}
  \put(8,60){\vector(1,0){10}}
\put(20,60){$\alpha_0$}
  \put(28,60){\vector(1,0){10}}
  \put(39,60){$\dots$}
  \put(48,60){\vector(1,0){10}}
  \put(60,60){$\alpha_0$}

\put(68,60){\vector(1,0){10}}
%\put(68,62){\vector(1,2){10}}
\put(68,58){\vector(1,-2){10}}

\put(80,59){$\alpha_2$}
\put(67,35){$\alpha_{p^s-1}$}
\put(105,75){\vector(1,0){10}}

\put(116,75){$\dots$}
\put(140,75){\vector(1,0){8}}
\put(125,75){\vector(1,0){8}}
\put(132,74){$\alpha_0$}
\put(148,74){$\alpha_1$}
\put(154,74){\vector(1,-1){12}}
\put(155,78){\vector(1,1){12}}
\put(155,76){\vector(1,0){12}}
\put(168,76){$\alpha_1$}
\put(168,88){$\alpha_0$}
\put(168,60){$\alpha_{p^s-1}$}

\put(85,63){\vector(1,1){12}}
\put(87,61){\vector(1,0){10}}
\put(85,59){\vector(1,-1){12}}
\put(98,74){$\alpha_0$}
\put(98,60){$\alpha_1$}
\put(98,45){$\alpha_{p^s-1}$}

\put(81,36){\vector(1,0){14}}
\put(83,34){\vector(1,-1){12}}
\put(81,32){\vector(0,-1){10}}
\put(96,35){$\alpha_0$}
\put(96,20){$\alpha_1$}
\put(78,18){$\alpha_{p^s-1}$}
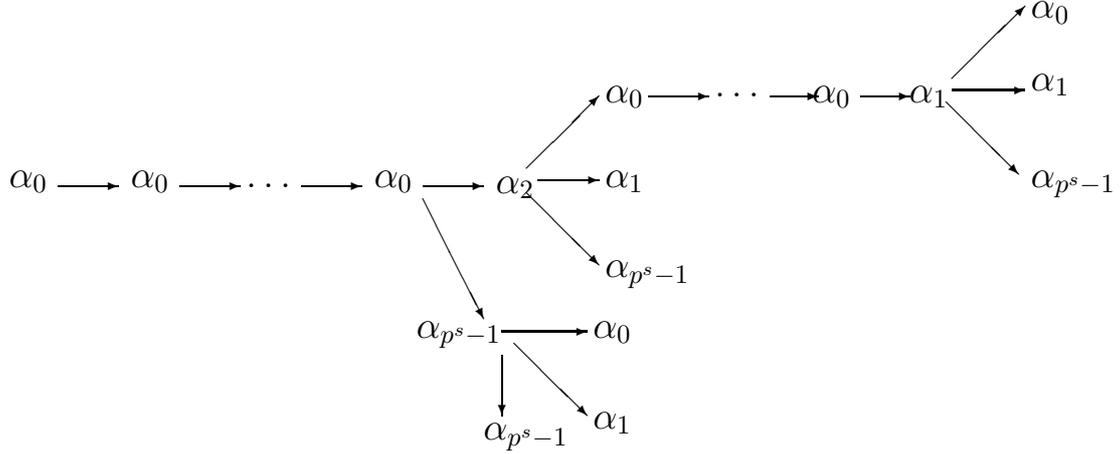
\put(40,5){Figure 5. Tree after moving }
      \end{picture}

 To obtain another $N$-valid trees we introduce the concept of {\it basic step} in the following way.\\
   1.Let $T$ be a $N$-valid tree.  Take  a subtree $T_{j(N-1+\nu),N-1+\nu}$ with a node   $\alpha_{j(N-1+\nu)}^{(N-1+\nu)},\ \nu\ge 1$ of the level $N-1+\nu$ as  a root. (see figures 3 and 4)\\
   2.Take a path
   $$
   \alpha_{j(\nu)}^{(\nu)}=\alpha_{j(N-1+\nu-N+1)}^{(N-1+\nu-N+1)}\rightarrow
   \dots \rightarrow \alpha_{j(N-1+\nu-1)}^{(N-1+\nu-1)}\rightarrow
   \alpha_{j(N-1+\nu)}^{(N-1+\nu )}
   $$
   of length $N-1$ which ends   in this node .  \\
   3. Find a path of the length $N-2$  which ends in leaf $\alpha=\alpha_{j(N-1+\nu-1)}^{(N-1+\nu-1)}$  and which coincides with the path
   $$
   \alpha_{j(N-1+\nu-N+1)}^{(N-1+\nu-N+1)}\rightarrow
   \dots \rightarrow \alpha_{j(N-1+\nu-1)}^{(N-1+\nu-1)}
   $$
      4. Move the subtree $T_{j(N-1+\nu),N-1+\nu}$ to the leaf $\alpha=\alpha_{j(N-1+\nu-1)}^{(N-1+\nu-1)}$. See Fig.5

   If the original tree was $N$-valid, then after the employment of  the basic step  we obtain a $N$-valid tree again. Thus,  applying  the basic algorithm to the basic tree finite number of times, we will obtain different $N$-valid trees.

Let $T$ be a $N$-valid tree. Choose any path of length greater than $N-1$
$$
P={\bf 0}_1\rightarrow\dots\rightarrow{\bf 0}_N\rightarrow \bar{\alpha}_\nu \rightarrow \bar{\alpha}_{\nu-1}\rightarrow\dots
\rightarrow\bar{\alpha}_{-N+1}\rightarrow\bar{\alpha}_{-N}=
$$
$$
= \bar{\alpha}_{\nu+N} \rightarrow
\bar{\alpha}_{\nu+N-1} \rightarrow\dots\rightarrow
\bar{\alpha}_{\nu+1} \rightarrow
\bar{\alpha}_\nu \rightarrow \bar{\alpha}_{\nu-1}\rightarrow\dots
\rightarrow\bar{\alpha}_{-N+1}\rightarrow\bar{\alpha}_{-N}.
$$
Since the tree $T$ is $N$-valid it follows that  $\bar{\alpha}_{\nu}\ne 0$.
Let us construct cosets
$$
(K^+_{-N})^\bot {\bf r}_{-N}^{\bar{\alpha}_{-N}}{\bf r}_{-N+1}^{\bar{\alpha}_{-N+1}}\dots{\bf r}_{-1}^{\bar{\alpha}_{-1}}
{\bf r}_{0}^{\bar{\alpha}_{0}},
$$
$$
 (K^+_{-N})^\bot
{\bf r}_{-N}^{\bar{\alpha}_{-N+1}}{\bf r}_{-N+1}^{\bar{\alpha}_{-N+2}}\dots
{\bf r}_{-1}^{\bar{\alpha}_{0}}
{\bf r}_{0}^{\bar{\alpha}_{1}},
$$
$$\dots \dots \dots \dots \dots \dots
$$
\begin{equation}                                \label{eq5.3}
(K^+_{-N})^\bot
{\bf r}_{-N}^{\bar{\alpha}_{\nu}}{\bf r}_{-N+1}^{\bar{\alpha}_{\nu+1}}\dots
{\bf r}_{-1}^{\bar{\alpha}_{\nu+N-1}}
{\bf r}_{0}^{\bar{\alpha}_{\nu+N}}
\end{equation}
and denote the union of all such cosets as $\tilde E$.
\begin{lemma}$\tilde E$ is $(N,1)$-elementary set.
\end{lemma}
{\it Proof.} Since $T$ is $N$-valid tree it follows that for any coset
$$
(K^+_{-N})^\bot
{\bf r}_{-N}^{\bar{\alpha}_{\nu}}{\bf r}_{-N+1}^{\bar{\alpha}_{\nu+1}}\dots
{\bf r}_{-1}^{\bar{\alpha}_{\nu+N-1}}\subset (K^+_{0})^\bot
$$
there exists unique shift
$$
(K^+_{-N})^\bot
{\bf r}_{-N}^{\bar{\alpha}_{\nu}}{\bf r}_{-N+1}^{\bar{\alpha}_{\nu+1}}\dots
{\bf r}_{-1}^{\bar{\alpha}_{\nu+N-1}}
{\bf r}_{0}^{\bar{\alpha}_{\nu+N}}\subset (K^+_{1})^\bot .
$$
It means that $\tilde E$ is $(N,1)$-elementary set. $\square$
\begin{definition}
 We say that the set $\tilde E_X$ is a
periodic extension of $\tilde E$ if
$$
\tilde E_X=\bigcup
\limits_{s=1}^\infty\bigsqcup\limits_{\bar{\alpha}_1,\dots,\bar{\alpha}_s=0}^{p-1}\tilde
E {\bf r}_1^{\bar{\alpha}_1}{\bf r}_2^{\bar{\alpha}_2}\dots {\bf r}_{\nu}^{\bar{\alpha}_{\nu}}.
$$
If $\bigcap\limits_{n=0}^\infty \tilde E_X{\cal A}^n=E$
then we say that  $\tilde E$ generates this set $E$, and the $N$-valid tree generates $E$ also.
\end{definition}
\begin{lemma}
Let $T$ be a $N$-valid  tree of height $H$.  Then
$$
\prod\limits_{n=0}^\infty{\bf 1}_{\tilde E_X}(\chi{\cal
A}^{-n})=\prod\limits_{n=0}^{H-N+1}{\bf 1}_{\tilde E_X}(\chi{\cal
A}^{-n})
$$
if $\chi\in (K^+_{H-2N+2})^\bot$.
\end{lemma}
{\bf Proof.} Since $\tilde E_X\supset (K^+_{-N})^\bot$ and $(K_{l}^+)^\bot{\cal A}=(K_{l+1}^+)^\bot$
it follows that
$$
{\bf 1}_{\tilde
E_X}((K_{H-2N+2}^+)^\bot{\cal A}^{-H+N-2})={\bf 1}_{\tilde
E_X}((K^+_{-N})^\bot)=1
$$
and the lemma is proved. $\square$

\begin{lemma}
Let $T$ be a $N$-valid  tree of height $H$. Suppose the tree $T$
generates the set  $E\subset X$. Then $E$ is
an $(N,H-2N+1)$-elementary set.
\end{lemma}
{\bf Proof.} Let us denote
$$
m(\chi)={\bf 1}_{\tilde
E_X}(\chi),\;\;M(\chi)=\prod\limits_{n=0}^\infty m(\chi{\cal
A}^{-n}).
$$
First we note that $M(\chi)={\bf 1}_E(\chi)$. Indeed
$$
{\bf 1}_E(\chi)=1\Leftrightarrow\chi\in E\Leftrightarrow
\forall\,n,\;\chi{\cal A}^{-n}\in\tilde
E_X\Leftrightarrow\forall\,n,\;{\bf 1}_{\tilde E_X}(\chi{\cal
A}^{-n})=1\Leftrightarrow
$$
$$
\forall\,n,\;m(\chi{\cal
A}^{-n})=1\Leftrightarrow\prod\limits_{n=0}^\infty m(\chi{\cal
A}^{-n})=1\Leftrightarrow M(\chi)=1.
$$
Now we will prove, that ${\bf 1}_E(\chi)=0$ for
$\chi\in(K^+_{ H-2N+2})^\bot\setminus (K^+_{ H-2N+1})^\bot $.\\
 In another words we need to prove that
$$
{\bf 1}_E(K_{-N}^+)^\bot
{\bf r}_{-N}^{\bar{\alpha}_{-N}}{\bf r}_{-N+1}^{\bar{\alpha}_{-N+1}}\dots
{\bf r}_{H-2N+1}^{\bar{\alpha}_{H-2N+1}})=0
$$
for $\bar{\alpha}_{H-2N+1}\ne 0$.

By the definition of cosets \eqref{eq5.3},
 $m((K_{-N}^+)^\bot
{\bf r}_{-N}^{\bar{\alpha}_{-N}}{\bf r}_{-N+1}^{\bar{\alpha}_{-N+1}}\dots {\bf r}_0^{\bar{\alpha}_0})
\ne 0$ if and only if
the vector
$(\bar{\alpha}_0,\bar{\alpha}_1,...,\bar{\alpha}_{-N+1},{\bar{\alpha}}_{-N})$ is a path
$(\bar{\alpha}_0\rightarrow\bar{\alpha}_1\rightarrow...\rightarrow
\bar{\alpha}_{-N+1}\rightarrow\bar{\alpha}_{-N})$
 of
the tree $T$.

 Since $\tilde E_X$ is a periodic
extension of $\tilde E$ it follows that the function $m(\chi)={\bf
1}_{\tilde E_X}(\chi)$ is periodic with any period
${\bf r}_{1}^{\bar{\alpha}_{1}}{\bf r}_{2}^{\bar{\alpha}_{2}}\dots {\bf r}_{\nu}^{\bar{\alpha}_{\nu}}$,
$\nu\in\mathbb N$, i.e. $m(\chi
{\bf r}_{1}^{\bar{\alpha}_{1}}{\bf r}_{2}^{\bar{\alpha}_{2}}\dots
{\bf r}_{\nu}^{\bar{\alpha}_{\nu}})=m(\chi)$ when $\chi\in (K^+_1)^\bot$.
Using this fact we can write $M(\chi)$ for  $\chi\in
(K^+_{H-2N+2})^\bot\setminus (K^+_{H-2N+1})^\bot$ in the form
$$
M((K^+_{-N})^\bot\zeta)=M((K^+_{-N})^\bot
{\bf r}_{-N}^{\bar{\alpha}_{-N}}{\bf r}_{-N+1}^{\bar{\alpha}_{-N+1}}\dots
{\bf r}_{H-2N+1}^{\bar{\alpha}_{H-2N+1}})=
$$
$$
 =m((K^+_{-N})^\bot {\bf r}_{-N}^{\bar{\alpha}_{-N}}
 {\bf r}_{-N+1}^{\bar{\alpha}_{-N+1}}\dots {\bf r}_{0}^{\bar{\alpha}_{0}})
 m((K^+_{-N})^\bot {\bf r}_{-N}^{\bar{\alpha}_{-N+1}}{\bf r}_{-N+1}^{\bar{\alpha}_{-N+2}}\dots
 {\bf r}_{0}^{\bar{\alpha}_{1}})\dots
$$
\begin{equation}  \label{eq5.4}
 m((K^+_{-N})^\bot
{\bf r}_{-N}^{\bar{\alpha}_{H-3N+1}}{\bf r}_{-N+1}^{\bar{\alpha}_{H-3N+2}}\dots
{\bf r}_{-1}^{\bar{\alpha}_{H-2N}}
 {\bf r}_{0}^{\bar{\alpha}_{H-2N+1}})
 \end{equation}
$$
m((K^+_{-N})^\bot
{\bf r}_{-N}^{\bar{\alpha}_{H-3N+2}}{\bf r}_{-N+1}^{\bar{\alpha}_{H-3N+3}}
r_{-1}^{\bar{\alpha}_{H-2N+1}})
$$
$$
\dots\dots\dots\dots\dots\dots\dots\dots\dots
$$
$$
m((K^+_{-N})^\bot {\bf r}_{-N}^{\bar{\alpha}_{H-2N}}{\bf r}_{-N+1}^{\bar{\alpha}_{H-2N+1}})m((K^+_{-N})^\bot {\bf r}_{-N}^{\bar{\alpha}_{H-2N+1}}).
$$
Assume that $M((K^+_{-N})^\bot\zeta)\neq 0$. Then all
factors in \eqref{eq5.4} are nonzero. So we have the path
$$
{\bf 0}\rightarrow\dots\rightarrow {\bf 0}\rightarrow\bar{\alpha} _{H-2N+1}\neq {\bf 0}\rightarrow \bar{\alpha}_{H-2N}\rightarrow\dots\rightarrow \bar{\alpha}
_0\rightarrow \dots \rightarrow\bar{\alpha}_{-N+1} \rightarrow\bar{\alpha}
_{-N},
$$
where there are $N$ zeroes at the beginning of the path. The length of such path is $H+1$, which contradicts the condition height of $T$ equals $H$.

Now we prove that $E$ is $(N,H-2N+1)$ elementary set.
Since the tree $T$ is $N$-valid, it has all possible combinations of $N$ elements  $\bar{\alpha}_i\in GF(p^s)$ as its paths, and we have the first property of elementary sets satisfied. Also, since height of $T$ is $H$, there exists a path
$$
\bar{\alpha}_1={\bf 0}\rightarrow\dots\rightarrow \bar{\alpha}_N={\bf 0}\rightarrow\bar{\alpha}_{N+1}\neq {\bf 0}\rightarrow\bar{\alpha}_{N+2}\rightarrow\dots\rightarrow\bar{\alpha}_{H+1}
$$
of length $H$. Such path generates cosets
$$
(K^+_{-N})^\bot {\bf r}_{-N}^{\bar{\alpha}_{N+l}}{\bf r}_{-N+1}^{\bar{\alpha}_{N+l-1}}\dots {\bf r}_{-N+l-1}^{\bar{\alpha}_{N+1}}
\subset (K^+_{-N+l})^\bot\setminus (K^+_{-N+l-1})^\bot
$$
for all $l=1,2,\dots,H+1-N$.
 Thus we can conclude that $E$ is $(N,H-2N+1)$-elementary set and the lemma is proved.
  $\square$.

 Now we can formulate an algorithm for constructing the refinable function that generates Riesz MRA.

 {\bf RF-algorithm.}\\
 \noindent
 1.Construct $N$-valid tree $T$ of height $H\ge 1$ using basic $N$-valid tree and basic steps.\\
 2.Construct the set $\tilde E \subset (K^+_1)^\bot$  using formulas (\ref{eq5.3}).\\
 3.Construct the function $m_0(\chi)$ on the set $ (K^+_1)^\bot$ such that

 3.1.$m_0((K^+_N)^\bot)=1$,

 3.2.${\rm supp}\  m_0(\chi)=\tilde E$,

 3.3.$0<A\le|m_0(\chi)|^2\le B$,\\
 4.  Extend the function $m_0$ periodically with any period ${\bf r}_1^{\bar{\alpha}_1}{\bf r}_2^{\bar{\alpha}_2}\dots {\bf r}_{\nu}^{\bar{\alpha}_{\nu}}$. It is evident that ${\rm supp}\ m_0(\chi)={\bf 1}_{\tilde
E_X}(\chi)$ and $A\le|m_0(\chi)|^2\le B$.\\
5. Set $\hat{\varphi}(\chi)=\prod_{n=0}^{H-N+1}m_0(\chi {\cal A}^{-n})$

{\bf Remark.} The Fourier transform $\hat{\varphi}$ may be calculated in the following way.
 Take any path
 $$
\bar{\alpha}_{\nu+N}\rightarrow\bar{\alpha}_{\nu+N-1}\rightarrow \dots \rightarrow\bar{\alpha}_{\nu+1}\rightarrow\bar{\alpha}_{\nu}\rightarrow\bar{\alpha}_{\nu-1}
\rightarrow\dots\rightarrow\bar{\alpha}_{-N}
$$
of the length $\ge N$ in that
$\bar{\alpha}_{\nu+N}=\bar{\alpha}_{\nu+N-1}=\dots =\bar{\alpha}_{\nu+1}=0,\bar{\alpha}_{\nu}\ne 0$. Then we set
$$
\hat{\varphi}((K^+_{-N})^\bot {\bf r}_{-N}^{\bar{\alpha}_{-N}}\dots{\bf r}_{\nu}^{\bar{\alpha}_{\nu}}\dots {\bf r}_{\nu+N}^{\bar{\alpha}_{\nu+N}})=
$$
$$
=m_0({\bf r}_{-N}^{\bar{\alpha}_{-N}}\dots{\bf r}_{0}^{\bar{\alpha}_{0}})
m_0({\bf r}_{-N}^{\bar{\alpha}_{-N+1}}\dots{\bf r}_{0}^{\bar{\alpha}_{1}})\dots
m_0({\bf r}_{-N}^{\bar{\alpha}_{\nu}}\dots{\bf r}_{0}^{\bar{\alpha}_{\nu+N}})
$$

 \begin{theorem}The function $\varphi$ generates Riesz MRA with constants $A^{H-N+2}$ and $ B^{H-N+2}$.
  \end{theorem}
 {\bf Proof.} It is evident that $A^{H-N+2}\le|\hat{\varphi}(\chi)|^2\le B^{H-N+2}$ and $\hat{\varphi}((K^+_{-N})^\bot)=1$. By lemma 5.3
 $$
 {\rm supp}\prod_{n=0}^{\infty}m_0(\chi {\cal A}^{-n})=E=\bigcap\limits_{n=0}^\infty \tilde E_X{\cal A}^n\subset (K^+_{H-2N+1})^\bot,
 $$
 so that
 $$
 {\rm supp}\prod_{n=0}^{\infty}m_0(\chi {\cal A}^{-n})=0
 $$
 on the set $(K^+_{H-2N+2})^\bot\setminus (K^+_{H-2N+1})^\bot$.
  Consequently, by lemma 5.2
  $$
 {\rm supp}\prod_{n=0}^{\infty}m_0(\chi {\cal A}^{-n})=
 {\rm supp}\prod_{n=0}^{H-N+1}m_0(\chi {\cal A}^{-n})=\hat{\varphi}(\chi).
 $$
 So, by theorem 4.2 the function $\varphi$ generates Riesz MRA. $\square$
\section{Construction of Riesz wavelets}
In this section we will give an algorithm for constructing wavelets.
We will use the result of B. Behera and  Q. Jahan \cite{BJ3}, which we formulate in our notations.\\

Let $\{V_j\}$ and $\{\tilde V_j\}$ be biorthogonal MRAs with scaling functions $\varphi,\ \tilde \varphi$ and masks $m_0(\chi), \tilde m_0(\chi)$ respectively.
Assume that there exist  periodic functions
$m_{\bf l}$ and $\tilde m_{\bf l},\  ( {\bf l}\in GF(p^s),\ {\bf l\neq {\bf 0}})$, such that for any ${\bf a}_{-N}\dots{\bf a}_{-1}\in GF(p^s)$ and for any $\chi\in (F_{-N}^{(s)})^\bot$
\begin{equation} \label{eq6.1}
\sum\limits_{{\bf a}_0\in GF(p^s)} m_{\bf k}(\chi {\bf r}_{-N}^{{\bf a}_{-N}}\dots {\bf r}_0^{{\bf a}_0})
 \overline{m_{\bf l}(\chi {\bf r}_{-N}^{{\bf a}_{-N}}\dots {\bf r}_0^{{\bf a}_0})}=
\delta_{{\bf k},{\bf l}}.
\end{equation}
 Define wavelets $\psi^{(\bf l)}$ and $\tilde\psi^{(\bf l)}$ by the equations
\begin{equation} \label{eq6.2}
\hat\psi^{(\bf l)}(\chi)=m_{\bf l}(\chi)\hat\varphi(\chi{\cal A}^{-1}),\
 \hat{\tilde\psi}^{(\bf l)}(\chi)=\tilde m_{\bf l}(\chi)\hat{\tilde\varphi}(\chi{\cal A}^{-1}),
\end{equation}

\begin{theorem}[\cite{BJ3}] Let $\varphi$ and $\tilde \varphi$ be the scaling functions for dual MRAs and $\psi_{\bf l}, \tilde\psi_{\bf l},  {\bf l}\in GF(p^s)$  be the associated wavelets satisfying the matrix condition (6.1). Then the collections
$$
\{\psi^{({\bf l})}_{n,h}=p^{\frac{ns}{2}}\psi^{({\bf l})}({\cal A}^n x\dot-h):h\in H_0,n\in \mathbb Z\}
$$
and
$$
\{\tilde{\psi}^{({\bf l})}_{n,h}=p^{\frac{ns}{2}}\tilde{\psi}^{({\bf l})}({\cal A}^n x\dot-h):h\in H_0,n\in \mathbb Z\}
$$ are biorthogonal. If addition
$$
|\hat\varphi((K_n^+)^\bot\setminus (K_{n-1}^+)^\bot )|\le\frac{C}{(1+p^{ns})^{\frac12+\varepsilon}},\
|\hat{\tilde\varphi}((K_n^+)^\bot\setminus (K_{n-1}^+)^\bot )|\le\frac{C}{(1+p^{ns})^{\frac12+\varepsilon}},
$$
$$
|\hat\psi^{(\bf l)}((K_n^+)^\bot)|\le Cp^{ns},\ |\hat{\tilde\psi}^{(\bf l)}((K_n^+)^\bot)|\le Cp^{ns},
$$
for some constant $C>0,\varepsilon >0$, then  systems $\{{\psi}^{({\bf l})}_{n,h}\}$ and $\{\tilde{\psi}^{({\bf l})}_{n,h}\}$ form Riesz bases for $L_2(K)$.
\end{theorem}

Now we can continue to  construct wavelets.
Let $T$ be $N$-valid tree. Using  RF-algorithm we construct functions $m_0(\chi)$, $\hat\varphi(\chi)$ and set
$$
\tilde m_0(\chi)=\left\{\begin{array}{ll}
                            0, & m_0(\chi)=0 \\
                            \frac{1}{m_0(\chi)}, & m_0(\chi)\ne 0
                          \end{array}
                        \right. ,\qquad
\hat{\tilde\varphi}(\chi)=\prod\limits_{n=0}^\infty\tilde m_0(\chi{\cal A}^{-n}).
$$
It is evident  $\hat\varphi(\chi)\overline{\hat{\tilde\varphi}(\chi)}=1$. Define functions
$$
m_{{\bf l}}(\chi)=m_0(\chi{\bf r}_0^{-{\bf l}}),\;\tilde m_{{\bf l}}(\chi)=\tilde m_0(\chi{\bf r}_0^{-{\bf l}}).
$$
\begin{lemma}
The following properties  are true\\
1) $m_{{\bf l}}(\tilde E_X{\bf r}_0^{{\bf l}})\ne 0, \tilde m_{{\bf l}}(\tilde E_X{\bf r}_0^{{\bf l}})\ne 0$ for any  ${\bf l}\in GF(p^s)$.\\
2) $m_{{\bf l}}(\tilde E_X{\bf r}_0^{{\bf a}})=\tilde m_{{\bf l}}(\tilde E_X{\bf r}_0^{{\bf a}})=0$ for  ${\bf l}\ne {\bf a}$.\\
3) $m_{\bf l}(E)=\tilde m_{\bf l}(E)=0$ for  ${\bf l}\ne {\bf 0}$.\\
4) $m_{\bf l}(\chi) m_{\bf k}(\chi)=\tilde m_{\bf l}(\chi)\tilde m_{\bf k}(\chi)=0$ for ${\bf k}\ne {\bf l}$.
\end{lemma}
{\bf Proof}. 1) If ${\bf l}\ne {\bf 0}$ then $m_{{\bf l}}(\tilde E_X{\bf r}_0^{{\bf l}})=m_0(\tilde E_X{\bf r}_0^{-{\bf l}}{\bf r}_0^{\bf l})=m_0(\tilde E_X)\ne 0$.\\
2) Let $\chi\in \tilde E$ and
 $\chi=\chi_{-N}{\bf r}_{-N}^{{\bf a}_{-N}}\dots{\bf r}_{-1}^{{\bf a}_{-1}}{\bf r}_{0}^{{\bf a}_{0}},\ \chi_{-N}\in (K_{-N}^+)^\bot$. It means that $m_0(\chi_{-N}{\bf r}_{-N}^{{\bf a}_{-N}}\dots{\bf r}_{-1}^{{\bf a}_{-1}}{\bf r}_{0}^{{\bf a}_{0}})\ne 0.$ Therefore
 if ${\bf a}\ne {\bf l}$, then
$$
m_{{\bf l}}((K_{-N}^+)^\bot {\bf r}_{-N}^{{\bf a}_{-N}}\dots{\bf r}_{-1}^{{\bf a}_{-1}}{\bf r}_{0}^{{\bf a}_{0}}{\bf r}_{0}^{{\bf a}})=m_0((K_{-N}^+)^\bot{\bf r}_{-N}^{{\bf a}_{-N}}\dots{\bf r}_{-1}^{{\bf a}_{-1}}{\bf r}_{0}^{{\bf a}_{0}-{\bf l}+{\bf a}}=0,
$$
since there cannot be two different paths from the node  ${\bf a}_{-N}$ to the root.\\
3) It follows from property  2) that  $m_{{\bf l}}(\tilde E_X)=m_{{\bf l}}(\tilde E_X{\bf r}_0^{{\bf 0}})=0$ for ${\bf l}\ne 0$.\\
4) If $\chi\in \tilde E$  then
$$
m_{{\bf l}}(\chi)=m_0(\chi_{-N}{\bf r}_{-N}^{{\bf\alpha}_{-N}}\dots{\bf r}_{-1}^{{\bf\alpha}_{-1}}{\bf r}_{0}^{{\bf\alpha}_{0}-{\bf l}}),
$$
$$
m_{{\bf k}}(\chi)=m_0(\chi_{-N}{\bf r}_{-N}^{{\bf\alpha}_{-N}}\dots{\bf r}_{-1}^{{\bf\alpha}_{-1}}{\bf r}_{0}^{{\bf\alpha}_{0}-{\bf k}}),
$$
 where $\chi_{-N}\in (K_{-N}^+)^\bot$.
 Since there cannot be two different paths from the node  ${\bf a}_{-N}$ to the root we see that  $m_{{\bf l}}(\chi)=0$, or $m_{{\bf k}}(\chi)=0$. $\square$

Define functions  $\psi^{({\bf l})}$ and $\tilde\psi^{({\bf l})}$ by equations  (6.2).
\begin{theorem}
1) The collectins $(\psi^{({\bf l})}_{n,h})$ and $(\tilde\psi^{({\bf l})}_{n,h})$ are biorthogonal.\\
2) The systems  $(\psi^{({\bf l})}_{n,h})$ and $(\tilde\psi^{({\bf l})}_{n,h})$ form Riesz bases for $L_2(K)$.
\end{theorem}
{\bf Proof.}
Check equality \eqref{eq6.1}. By Lemma 6.1, so that
 $m_{\bf l}m_{\bf k}=0$ when  ${\bf k}\ne {\bf l}$. Therefore, it suffices to prove the equation
\begin{equation} \label{eq6.3}
\sum_{\alpha_0\in GF(p^s)}m_{\bf l}(\chi_{-N}{\bf r}_{-N}^{{\bf\alpha}_{-N}}\dots {\bf r}_{-1}^{{\bf\alpha}_{-1}}{\bf r}_{0}^{{\bf\alpha}_{0}})\overline{\tilde m_{\bf l}(\chi_{-N}{\bf r}_{-N}^{{\bf\alpha}_{-N}}\dots {\bf r}_{-1}^{{\bf\alpha}_{-1}}{\bf r}_{0}^{{\bf\alpha}_{0}})}=1.
\end{equation}
Since there cannot be two different paths from the node  ${\bf a}_{-N}$ to the root, it follows that
   (6.3)    includes only one non-zero term is equal to one.
    By lemma  6.1 $m_{\bf l}(E)=0$ for  ${\bf l}\ne 0$. Therefore $\psi^{({\bf l})}((K_{-N}^+)^\bot)=0$.

Since
${\rm supp}\, \hat\varphi(\chi)\subset (K_{H-2N+1}^+)^\bot$, it follows that  ${\rm supp}\, \hat\psi^{({\bf l})}(\chi)\subset (K_{H-2N+2}^+)^\bot$. By analogy, $\tilde \psi^{({\bf l})}((K_{-N}^+)^\bot)=0$ and ${\rm supp}\, \hat{\tilde{\psi}}^{({\bf l})}(\chi)\subset (K_{H-2N+2}^+)^\bot$.
So all conditions of theorem  6.1 is fulfilled, and theorem  6.2 is proved. $\square$

 Finally we can write an algorithm to construct Riesz-wavelets.\\
{\bf W-algorithm.}\\
1) Construct $N$-valid tree $T$ using the basic steps.\\
2) Construct the mask  $m_0(\chi)$ and refinable function $\varphi(\chi)$  using RF-algorithm.\\
3) Define functions  $m_{{\bf l}}(\chi)=m_0(\chi {\bf r}_0^{-{\bf l}})$.\\
4) Set $\hat\psi^{({\bf l})}(\chi)=m_{{\bf l}}(\chi)\hat\varphi(\chi{\cal A}^{-1})$.\\
5)Find wavelets $\psi^{({\bf l})}(\chi)$ using inverse Fourier transform.

This research was carried out with the financial
support the Russian Foundation for Basic Research (grant
no.~16-01-00152)

 \end{document}